\documentclass{article}
\usepackage{amssymb,amsmath,amsthm,graphicx}
\usepackage[all,color]{xy}

\def\qed{\hfill {\hbox{${\vcenter{\vbox{               
   \hrule height 0.4pt\hbox{\vrule width 0.4pt height 6pt
   \kern5pt\vrule width 0.4pt}\hrule height 0.4pt}}}$}}}

\def\utr{\, \underline{\triangleright}\, }
\def\otr{\, \overline{\triangleright}\, }

\newtheorem{theorem}{Theorem}

\newtheorem{corollary}[theorem]{Corollary}

\theoremstyle{definition}
\newtheorem{example}{Example}
\newtheorem{definition}{Definition}
\newtheorem{remark}{Remark}

\date{}

\title{\Large \textbf{Biquandle Module Invariants of Oriented Surface-Links}}

\author{Yewon Joung\footnote{Email: yjoung@msu.edu.}
\and Sam Nelson\footnote{Email: Sam.Nelson@cmc.edu. Partially supported by Simons Foundation collaboration grant 316709}}

\begin{document}
\maketitle

\begin{abstract}
We define invariants of oriented surface-links by enhancing the biquandle
counting invariant using \textit{biquandle modules}, algebraic structures
defined in terms of biquandle actions on commutative rings analogous to
Alexander biquandles. We show that bead colorings of marked graph diagrams
are preserved by Yoshikawa moves and hence define enhancements of the biquandle
counting invariant for surface links. We provide examples illustrating the
computation of the invariant and demonstrate that these invariants are not 
determined by the first and second Alexander elementary ideals and 
characteristic polynomials.
\end{abstract}

\parbox{5.5in} {\textsc{Keywords:} Biquandle modules, counting invariants, 
surface-links, marked graph diagrams

\smallskip

\textsc{2010 MSC:} 57M27, 57M25}

\section{Introduction}\label{I}

In \cite{AG} the notion of \textit{quandle modules} was introduced and used
in \cite{CEGS} 
to generalize quandle cocycle invariants of oriented classical knots
and links using \textit{dynamical cocycles}. In later
papers \cite{BN,BKLNS,CN,GN,HHNYZ,NP} the second listed 
author and collaborators adapted the dynamical cocycle idea to various 
settings in classical, virtual and twisted virtual knot theory,
defining enhancements of the quandle, rack and biquandle counting invariants.
Quandle module colorings of knots can be understood as secondary colorings
of quandle colored knots using ``beads'' which obey a kind of
customized Alexander quandle coloring rule with coefficients depending on the
base quandle coloring. In particular, quandle modules can be understood as
generalized Alexander quandles for quandle-colored knots and links.

In \cite{KJL} the first listed author and coauthors considered Alexander
biquandle colorings of oriented surface-links represented by \textit{marked
graph diagrams}, also known as \textit{marked vertex diagrams} or 
\textit{ch-diagrams}. In particular, the methods of \cite{KJL} distinguished 
most of the oriented surface-links of small $ch$-index as identified in 
\cite{Y} but did not distinguish the surface-links $6_1^{0,1}$ and $8_1^{1,1}$. 

In this paper we consider biquandle module invariants in the setting of
orientable surface-links. We show that biquandle module colorings are preserved
by a generating set of Yoshikawa moves, and hence define invariants of oriented
surface-links. As an application, we exhibit a biquandle
module invariant which does distinguish $6_1^{0,1}$ and $8_1^{1,1}$,
showing in particular that biquandle module invariants are not determined
by the first and second Alexander elementary ideals and characteristic
polynomials. We provides explicit examples to illustrate the computation of 
the invariant and report the results of computer calculation of the invariant
for some choices of biquandle modules for the oriented links of small 
$ch$-index.

The paper is organized as follows. In Section \ref{M} we review the basics
of marked graph diagrams and surface-link theory. In Section \ref{BR}
we review biquandles and biquandle colorings of marked graph diagrams. 
In Section \ref{C} we recall the definition of biquandle modules (updated 
with current notation) and show that biquandle module ``bead'' colorings
are preserved by oriented Yoshikawa moves. We define biquandle module 
enhancement invariants for oriented surface-links and compute some 
examples. We end in Section \ref{Q} with some questions for future work.

\section{Marked Graph Diagrams}\label{M}

In this section, we review (oriented) marked graph diagrams representing 
surface-links.

\begin{definition}
A \textit{marked graph diagram}, also called a \textit{marked vertex diagram} or
\textit{ch-diagram}, is a planar 4-regular
graph diagram with vertices decorated as \textit{classical crossings}
and \textit{saddle crossings} as depicted. 
\[\includegraphics{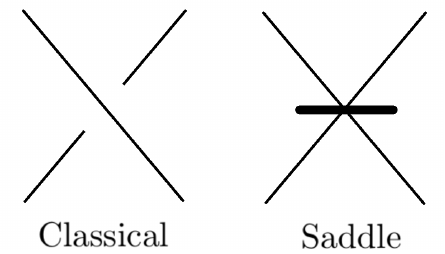}\]
A marked graph diagram is
\textit{orientable} if each edge in the graph can be directed such that
the classical crossings receive ``pass-through'' orientations and the saddle
crossings receive ``source-sink'' orientations. 
\[\includegraphics{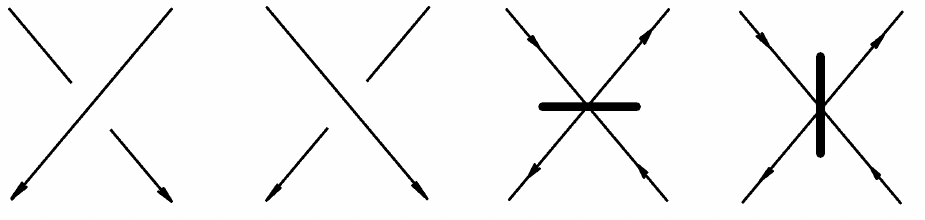}\]
A marked graph diagram
is \textit{admissible} if the two diagrams resulting from smoothing all saddle
points with the bars and against the bars are unlinks.
\end{definition}

A marked graph diagram represents an orientable surface-link in 
$\mathbb{R}^4$ if the diagram is admissible; the surface is obtained by
replacing each saddle crossing with a saddle as depicted and capping off 
the resulting unlinks above and below. 
\[\includegraphics{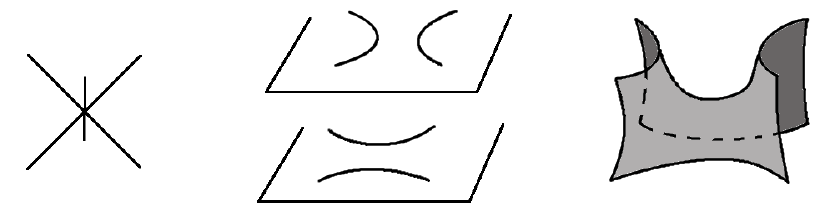}\]
Non-admissible marked graph diagrams
represent cobordisms between the links obtained through smoothing. 
Non-orientable admissible diagrams represent non-orientable surface-links,
and non-admissible non-orientable diagrams represent non-orientable
cobordisms. In particular, we may identify a classical link $L$ with the
surface-link given by the trivial cobordism (i.e., $L\times[0,1]$).

Two marked graph diagrams represent ambient isotopic surface-links in 
$\mathbb{R}^4$ if and only if they are related by the \textit{Yoshikawa moves}.
In \cite{KJL2} the first author and coauthors identified the generating set of
oriented Yoshikawa moves pictured here. 

\centerline{ 
\xy (12,2);(16,6) **@{-}, 
(12,6);(13.5,4.5) **@{-},
(14.5,3.5);(16,2) **@{-}, 
(16,6);(22,6) **\crv{(18,8)&(20,8)},
(16,2);(22,2) **\crv{(18,0)&(20,0)}, (22,6);(22,2) **\crv{(23.5,4)},
(7,8);(12,6) **\crv{(10,8)}, (7,0);(12,2) **\crv{(10,0)}, (12.4,5) *{\ulcorner}, (15.6,5) *{\urcorner},
(35,5);(45,5) **@{-} ?>*\dir{>}, (35,3);(45,3) **@{-} ?<*\dir{<},
(63.8,2) *{\urcorner},
(57,8);(57,0) **\crv{(67,7)&(67,1)}, (-5,4)*{\Gamma_1 :}, (73,4)*{},
\endxy }

\vskip.2cm

\centerline{ 
\xy (12,2);(16,6) **@{-}, 
(12,6);(13.5,4.5) **@{-},
(14.5,3.5);(16,2) **@{-}, 
(16,6);(22,6) **\crv{(18,8)&(20,8)},
(16,2);(22,2) **\crv{(18,0)&(20,0)}, (22,6);(22,2) **\crv{(23.5,4)},
(7,8);(12,6) **\crv{(10,8)}, (7,0);(12,2) **\crv{(10,0)}, (12.4,2.5) *{\llcorner}, (15.6,2.5) *{\lrcorner},
(35,5);(45,5) **@{-} ?>*\dir{>}, (35,3);(45,3) **@{-} ?<*\dir{<},(63.85,5.2) *{\lrcorner},
(57,8);(57,0) **\crv{(67,7)&(67,1)}, (-5,4)*{\Gamma_1' :}, (73,4)*{},
\endxy }

\vskip.2cm

\centerline{ \xy (7,7);(7,1)  **\crv{(23,6)&(23,2)}, (16,6.3);(23,7)
**\crv{(19,6.9)}, (16,1.7);(23,1) **\crv{(19,1.1)},
(14,5.7);(14,2.3) **\crv{(8,4)}, (10,6.9) *{<}, (20,6.9) *{>},
(35,5);(45,5) **@{-} ?>*\dir{>}, (35,3);(45,3) **@{-} ?<*\dir{<},
(57,7);(57,1) **\crv{(65,6)&(65,2)}, (73,7);(73,1)
**\crv{(65,6)&(65,2)}, (59,6.7) *{<}, (71,6.7) *{>}, (-5,4)*{\Gamma_2 :},
\endxy}

\vskip.2cm

\centerline{ 
\xy (7,6);(23,6) **\crv{(15,-2)}, 
(10,0);(11.5,1.8) **@{-}, 
(17.5,3);(14.5,6.6) **@{-},
(14.5,6.6);(10,12) **@{-}, 
(20,12);(15.5,6.6) **@{-},
(14.5,5.5);(12.5,3) **@{-},
(18.5,1.8);(20,0) **@{-},
(19.5,11) *{\urcorner}, 
(19,1.1) *{\lrcorner},  
(9,3.5) *{\ulcorner},
(35,7);(45,7) **@{-} ?>*\dir{>}, 
(35,5);(45,5) **@{-} ?<*\dir{<},
(57,6);(73,6) **\crv{(65,14)}, 
(70,12);(68.5,10.2) **@{-}, 
(67.5,9);(65.5,6.5) **@{-}, 
(64.6,5.5);(60,0) **@{-}, 
(62.5,9);(64.4,6.6) **@{-}, 
(64.4,6.6);(70,0) **@{-}, 
(61.5,10.2);(60,12) **@{-},
(69.5,11) *{\urcorner}, 
(69,1.1) *{\lrcorner},  
(58,7) *{\llcorner},
(-5,6)*{\Gamma_3:},
\endxy}

\vskip.2cm

 \centerline{ \xy 
 (7,6);(23,6)  **\crv{(15,-2)}, 
 (10,0);(11.5,1.8) **@{-},
(12.5,3);(20,12) **@{-}, 
(10,12);(17.5,3) **@{-}, 
(18.5,1.8);(20,0) **@{-}, 
(13,6);(17,6) **@{-}, (13,6.1);(17,6.1) **@{-}, (13,5.9);(17,5.9)
**@{-}, (13,6.2);(17,6.2) **@{-}, (13,5.8);(17,5.8) **@{-},
  (13,3.1) *{\urcorner}, (17.5,8.9) *{\llcorner}, (19,1.1) *{\lrcorner},  (21,3.5) *{\urcorner},(13,8) *{\ulcorner},
(35,7);(45,7) **@{-} ?>*\dir{>}, 
(35,5);(45,5) **@{-} ?<*\dir{<},
(57,6);(73,6)  **\crv{(65,14)}, 
(70,12);(68.5,10.2) **@{-},
(67.5,9);(60,0) **@{-}, 
(70,0);(62.5,9) **@{-}, 
(61.5,10.2);(60,12) **@{-}, 
(63,6);(67,6) **@{-}, (63,6.1);(67,6.1) **@{-}, (63,5.9);(67,5.9)
**@{-}, (63,6.2);(67,6.2) **@{-}, (63,5.8);(67,5.8) **@{-},
(62,2) *{\urcorner}, (67,8.3) *{\llcorner}, (67.5,3) *{\lrcorner},  (70.9,8) *{\lrcorner},(61.3,10) *{\ulcorner}, 
(-5,6)*{\Gamma_4:},
\endxy}

\vskip.2cm

 \centerline{ \xy 
  (13,2.2);(17,2.2)  **\crv{(15,1.7)}, 
  (7,6);(11,3)  **\crv{(10,3.5)}, 
  (23,6);(19,3)  **\crv{(20,3.5)}, 
 (10,0);(20,12) **@{-}, 
(10,12);(20,0) **@{-}, 
(13,6);(17,6) **@{-}, (13,6.1);(17,6.1) **@{-}, (13,5.9);(17,5.9)
**@{-}, (13,6.2);(17,6.2) **@{-}, (13,5.8);(17,5.8) **@{-}, 
(13,3.1) *{\urcorner}, (17.5,8.9) *{\llcorner}, (19,1.1) *{\lrcorner},  (21,3.5) *{\urcorner},(13,8) *{\ulcorner},
(35,7);(45,7) **@{-} ?>*\dir{>}, 
(35,5);(45,5) **@{-} ?<*\dir{<},
   (63,9.8);(67,9.8)  **\crv{(65,10.3)}, 
  (57,6);(61,9)  **\crv{(60,8.5)}, 
  (73,6);(69,9)  **\crv{(70,8.5)}, 
(70,12);(60,0) **@{-}, 
(70,0);(60,12) **@{-}, 
(63,6);(67,6) **@{-}, (63,6.1);(67,6.1) **@{-}, (63,5.9);(67,5.9)
**@{-},(63,6.2);(67,6.2) **@{-}, (63,5.8);(67,5.8) **@{-}, 
(62,2) *{\urcorner}, (67,8.3) *{\llcorner}, (67.5,3) *{\lrcorner},  (70.9,8) *{\lrcorner},(61.3,10) *{\ulcorner}, 
(-5,6)*{\Gamma_4':},
\endxy}

\vskip.2cm

 \centerline{ \xy (9,2);(13,6) **@{-}, (9,6);(10.5,4.5) **@{-},
(11.5,3.5);(13,2) **@{-}, (17,2);(21,6) **@{-}, (17,6);(21,2)
**@{-}, (13,6);(17,6) **\crv{(15,8)}, (13,2);(17,2) **\crv{(15,0)},
(7,7);(9,6) **\crv{(8,7)}, (7,1);(9,2) **\crv{(8,1)}, (23,7);(21,6)
**\crv{(22,7)}, (23,1);(21,2) **\crv{(22,1)}, 
(17,4);(21,4) **@{-}, (17,4.1);(21,4.1) **@{-}, (17,3.9);(21,3.9)
**@{-}, (17,4.2);(21,4.2) **@{-}, (17,3.8);(21,3.8) **@{-},
(10,3) *{\llcorner},  (12,3) *{\lrcorner}, (21,6) *{\llcorner},(21,2.2) *{\lrcorner},
(35,5);(45,5) **@{-} ?>*\dir{>}, (35,3);(45,3) **@{-} ?<*\dir{<},
(59,2);(63,6) **@{-}, (59,6);(63,2) **@{-}, (67,2);(71,6) **@{-},
(67,6);(68.5,4.5) **@{-}, (69.5,3.5);(71,2) **@{-}, (63,6);(67,6)
**\crv{(65,8)}, (63,2);(67,2) **\crv{(65,0)}, (57,7);(59,6)
**\crv{(58,7)}, (57,1);(59,2) **\crv{(58,1)}, (73,7);(71,6)
**\crv{(72,7)}, (73,1);(71,2) **\crv{(72,1)}, 
(63,4);(59,4) **@{-}, (63,4.1);(59,4.1) **@{-}, (63,3.9);(59,3.9)
**@{-}, (63,4.2);(59,4.2) **@{-}, (63,3.8);(59,3.8) **@{-},
(59.5,2.5) *{\llcorner},  (59,6) *{\lrcorner}, (71,6) *{\llcorner},(71,2.2) *{\lrcorner},
 (-5,4)*{\Gamma_5:},
\endxy }

\vskip.2cm

\centerline{ \xy (12,6);(16,2) **@{-}, (12,2);(16,6) **@{-},
(16,6);(22,6) **\crv{(18,8)&(20,8)}, (16,2);(22,2)
**\crv{(18,0)&(20,0)}, (22,6);(22,2) **\crv{(23.5,4)}, (7,8);(12,6)
**\crv{(10,8)}, (7,0);(12,2) **\crv{(10,0)}, (11,0.4) *{\urcorner}, (11,6) *{\ulcorner},(21.5,1)*{\urcorner},
(35,5);(45,5) **@{-} ?>*\dir{>}, (35,3);(45,3) **@{-} ?<*\dir{<},
(57,8);(57,0) **\crv{(67,7)&(67,1)}, (-5,4)*{\Gamma_6 :}, (73,4)*{},
(14,6);(14,2) **@{-}, (14.1,6);(14.1,2) **@{-}, (13.9,6);(13.9,2)
**@{-}, (14.2,6);(14.2,2) **@{-}, (13.8,6);(13.8,2) **@{-}, 
(63.8,2) *{\urcorner},
\endxy}

\vskip.2cm

\centerline{ \xy (12,6);(16,2) **@{-}, (12,2);(16,6) **@{-},
(16,6);(22,6) **\crv{(18,8)&(20,8)}, (16,2);(22,2)
**\crv{(18,0)&(20,0)}, (22,6);(22,2) **\crv{(23.5,4)}, (7,8);(12,6)
**\crv{(10,8)}, (7,0);(12,2) **\crv{(10,0)}, (11,0.4) *{\urcorner}, (11,6) *{\ulcorner},(21.5,1)*{\urcorner},
(35,5);(45,5) **@{-} ?>*\dir{>}, (35,3);(45,3) **@{-} ?<*\dir{<},
(57,8);(57,0) **\crv{(67,7)&(67,1)}, (-5,4)*{\Gamma'_6 :},
(73,4)*{}, (12,4);(16,4) **@{-}, (12,4.1);(16,4.1) **@{-},
(12,4.2);(16,4.2) **@{-}, (12,3.9);(16,3.9) **@{-},
(12,3.8);(16,3.8) **@{-}, (63.8,2) *{\urcorner},
\endxy}

\vskip.2cm

\centerline{ \xy (9,4);(17,12) **@{-}, (9,8);(13,4) **@{-},
(17,12);(21,16) **@{-}, (17,16);(21,12) **@{-}, (7,0);(9,4)
**\crv{(7,2)}, (7,12);(9,8) **\crv{(7,10)}, (15,0);(13,4)
**\crv{(15,2)}, (17,16);(15,20) **\crv{(15,18)}, (21,16);(23,20)
**\crv{(23,18)}, (21,12);(23,8) **\crv{(23,10)}, (7,12);(7,20)
**@{-}, (23,8);(23,0) **@{-},
(11,4);(11,8) **@{-}, 
(10.9,4);(10.9,8) **@{-}, 
(11.1,4);(11.1,8) **@{-}, 
(10.8,4);(10.8,8) **@{-}, 
(11.2,4);(11.2,8) **@{-},
(17,14);(21,14) **@{-}, 
(17,14.1);(21,14.1) **@{-},
(17,13.9);(21,13.9) **@{-}, 
(17,14.2);(21,14.2) **@{-},
(17,13.8);(21,13.8) **@{-},
(7,15) *{\wedge},(23,5) *{\wedge},(15,10) *{\llcorner},(8,2.3) *{\urcorner}, (21.5,16) *{\urcorner},(16,17) *{\lrcorner},(13.7,3) *{\lrcorner},
(35,11);(45,11) **@{-} ?>*\dir{>}, (35,9);(45,9) **@{-} ?<*\dir{<},
(71,4);(63,12) **@{-}, (71,8);(67,4) **@{-}, (63,12);(59,16) **@{-},
(63,16);(59,12) **@{-}, (73,0);(71,4) **\crv{(73,2)}, (73,12);(71,8)
**\crv{(73,10)}, (65,0);(67,4) **\crv{(65,2)}, (63,16);(65,20)
**\crv{(65,18)}, (59,16);(57,20) **\crv{(57,18)}, (59,12);(57,8)
**\crv{(57,10)}, (73,12);(73,20) **@{-}, (57,8);(57,0) **@{-},
(61,12);(61,16) **@{-}, 
(61.1,12);(61.1,16) **@{-},
(60.9,12);(60.9,16) **@{-}, 
(61.2,12);(61.2,16) **@{-},
(60.8,12);(60.8,16) **@{-},
(57,5) *{\wedge},(73,15) *{\wedge},(65,10) *{\lrcorner},(58,17) *{\ulcorner}, (71.5,3) *{\ulcorner},(66.3,3) *{\llcorner},(63.7,16.8) *{\llcorner},
(67,6);(71,6) **@{-}, 
(67,6.1);(71,6.1) **@{-}, 
(67,5.9);(71,5.9) **@{-}, 
(67,6.2);(71,6.2) **@{-},
(67,5.8);(71,5.8) **@{-},  
(-5,10)*{\Gamma_7:}, 
 \endxy}

\vskip.2cm

\centerline{ \xy (7,20);(14.2,11) **@{-}, (15.8,9);(17.4,7) **@{-},
(19,5);(23,0) **@{-}, (13,20);(7,12) **@{-}, (7,12);(11.2,7) **@{-},
(12.7,5.2);(14.4,3.2) **@{-}, (15.7,1.6);(17,0) **@{-},
(17,20);(23,12) **@{-}, (13,0);(23,12) **@{-}, (7,0);(23,20) **@{-},
(10,18);(10,14) **@{-}, (10.1,18);(10.1,14) **@{-},
(9.9,18);(9.9,14) **@{-}, (10.2,18);(10.2,14) **@{-},
(9.8,18);(9.8,14) **@{-}, (18,16);(22,16) **@{-},
(18,16.1);(22,16.1) **@{-}, (18,15.9);(22,15.9) **@{-},
(18,16.2);(22,16.2) **@{-}, (18,15.8);(22,15.8) **@{-},
 (8.5,17.7) *{\ulcorner}, (18.5,17.7) *{\ulcorner}, (9.1,9) *{\ulcorner}, (12,18.4) *{\llcorner}, (14.5,1.6) *{\llcorner}, (21.8,18.4) *{\llcorner}, (17,12) *{\urcorner}, (21.7,1.6) *{\lrcorner},
(35,11);(45,11) **@{-} ?>*\dir{>}, (35,9);(45,9) **@{-} ?<*\dir{<},
(73,20);(65.8,11) **@{-}, (64.2,9);(62.6,7) **@{-}, (61,5);(57,0)
**@{-}, (67,20);(73,12) **@{-}, (73,12);(68.8,7) **@{-},
(67.3,5.2);(65.6,3.2) **@{-}, (64.3,1.6);(63,0) **@{-},
(63,20);(57,12) **@{-}, (67,0);(57,12) **@{-}, (73,0);(57,20)
**@{-},
 (58.5,17.7) *{\ulcorner}, (68.5,17.7) *{\ulcorner}, (59.1,9) *{\ulcorner}, (62,18.4) *{\llcorner}, (64,1.2) *{\llcorner}, (71.8,18.4) *{\llcorner}, (67,12) *{\urcorner}, (71.7,1.6) *{\lrcorner},
(60,18);(60,14) **@{-}, (60.1,18);(60.1,14) **@{-},
(59.9,18);(59.9,14) **@{-}, (60.2,18);(60.2,14) **@{-},
(59.8,18);(59.8,14) **@{-}, (68,16);(72,16) **@{-},
(68,16.1);(72,16.1) **@{-}, (68,15.9);(72,15.9) **@{-},
(68,16.2);(72,16.2) **@{-}, (68,15.8);(72,15.8) **@{-},
(-5,10)*{\Gamma_8:}, 
\endxy}

Marked graph diagrams and Yoshikawa moves provide a very convenient calculus
for computing invariants of closed surface-links as well as cobordisms between
knots and links. In particular, a classical knot or link diagram $L$ considered
as a marked graph diagram can be pictured as a product $L\times[0,1]$. See 
\cite{KJL,KJL2} for more.

\section{Biquandles and  Colorings}\label{BR}

In this section we recall \textit{biquandles} and the biquandle counting 
invariant for marked graph diagrams.

\begin{definition}
Let $X$ be a set. A \textit{biquandle structure} on $X$ consists of two 
binary operations $\otr,\utr$ on $X$ satisfying for all $x,y,z\in X$
\begin{itemize}
\item[(i)] $x\utr x=x\otr x$,
\item[(ii)] the maps $\alpha_x,\beta_x:X\to X$ and $S:X\times X\to X\times X$
defined by $\alpha_x(y)=y\otr x$, $\beta_x(y)=y\utr x$ and
$S(x,y)=(y\otr x,x\utr y)$ are invertible, and
\item[(iii)] the \textit{exchange laws} are satisfied:
\[\begin{array}{rcl}
(x\utr y)\utr(z\utr y) & = & (x\utr z)\utr (y\otr z), \\
(x\utr y)\otr(z\utr y) & = & (x\otr z)\utr (y\otr z), \ \mathrm{and} \\
(x\otr y)\otr(z\otr y) & = & (x\otr z)\otr (y\utr z).
\end{array}\]
\end{itemize}
\end{definition}

\begin{example}
For any set $X$ and bijection $\sigma:X\to X$ the operations 
$x\utr y=x\otr y =\sigma(x)$ define a biquandle structure called a 
\textit{constant action biquandle}. If $\sigma$ is the identity then we
have a \textit{trivial biquandle}.
\end{example}

\begin{example}
Let $X$ be any module over $\mathbb{Z}[t^{\pm 1},s^{\pm 1}]$. Then
$X$ is a biquandle with operations
\[x\utr y=tx+(s-t)y,\quad x\otr y=sx\]
known as an \textit{Alexander biquandle}. Equivalently, any abelian group
$X$ with automorphisms $t,s:X\to X$ is an Alexander biquandle with 
\[x\utr y=t(x)+s(y)-t(y),\quad x\otr y=s(x).\]
\end{example}

\begin{example}
Given an oriented marked graph diagram $L$, let $\{x_1,\dots, x_n\}$ be 
a set of generators associated to the semiarcs in the diagram. A 
\textit{biquandle word} is either a generator or obtained recursively from the 
generators as $a\utr b$ or $a\otr b$ where $a,b$ are biquandle words.
Then the \textit{fundamental biquandle} of $L$, $\mathcal{B}(L)$, is the 
set of equivalence classes of biquandle words under the equivalence relation
generated by the biquandle axioms and the \textit{crossing relations}:
\end{example}

\[ \xy (10,-8);(26,8) **@{-} ?<*\dir{<},
(17,1);(10,8) **@{-}, (26,-8);(19,-1) **@{-} ?<*\dir{<},
(8,-10)*{y}, (29,-10)*{x \utr y}, (7.5,10)*{x}, (29,10)*{y \otr x},
(26,0)*{_{\varepsilon_i=+}}, 
(13,0)*{_{c_i}},
\endxy ~~~
\xy (10,8);(26,-8) **@{-} ?>*\dir{>}, 
(17,-1);(10,-8) **@{-} ?>*\dir{>}, 
(26,8);(19,1) **@{-},
(8,-10)*{x}, (29,-10)*{y \otr x}, (7.5,10)*{y}, (29,10)*{x \utr y},
(26,0)*{_{\varepsilon_i=-}}, 
(13,0)*{_{c_i}},
\endxy ~~~
\xy (10,8);(26,-8) **@{-} ?>*\dir{>}, 
(26,8);(10,-8) **@{-} ?>*\dir{>}, 
(8,-10)*{x}, (29,-10)*{x}, (7.5,10)*{x}, (29,10)*{x},
(13,0)*{_{v_i}},
\endxy\]

\centerline{ \xy 
 (7,6);(23,6)  **\crv{(15,-2)}, 
 (10,0);(11.5,1.8) **@{-},
(12.5,3);(20,12) **@{-}, 
(10,12);(17.5,3) **@{-}, 
(18.5,1.8);(20,0) **@{-}, 
(13,6);(17,6) **@{-}, (13,6.1);(17,6.1) **@{-}, (13,5.9);(17,5.9)
**@{-}, (13,6.2);(17,6.2) **@{-}, (13,5.8);(17,5.8) **@{-},
  (13,3.1) *{\urcorner}, (17.5,8.9) *{\llcorner}, (19,1.1) *{\lrcorner},  (21,3.5) *{\urcorner},(13,8) *{\ulcorner},
(9,14)*{x \utr y}, (5,8)*{y \otr x}, (8,0)*{x}, (15,0)*{y}, (22,0)*{x}, (25,8)*{y \otr x},
(35,7);(45,7) **@{-} ?>*\dir{>}, 
(35,5);(45,5) **@{-} ?<*\dir{<},
(57,6);(73,6)  **\crv{(65,14)}, 
(70,12);(68.5,10.2) **@{-},
(67.5,9);(60,0) **@{-}, 
(70,0);(62.5,9) **@{-}, 
(61.5,10.2);(60,12) **@{-}, 
(63,6);(67,6) **@{-}, (63,6.1);(67,6.1) **@{-}, (63,5.9);(67,5.9)
**@{-}, (63,6.2);(67,6.2) **@{-}, (63,5.8);(67,5.8) **@{-},
(62,2) *{\urcorner}, (67,8.3) *{\llcorner}, (67.5,3) *{\lrcorner},  (70.9,8) *{\lrcorner},(61.3,10) *{\ulcorner}, 
(60,14)*{x \utr y}, (56,5)*{y \otr x}, (65,3)*{x}, (65,12)*{y}, (77,5)*{y \otr x}, (70,14)*{x \utr y}, 
(-5,6)*{\Gamma_4:},
\endxy}

\vskip.3cm

 \centerline{ \xy 
  (13,2.2);(17,2.2)  **\crv{(15,1.7)}, 
  (7,6);(11,3)  **\crv{(10,3.5)}, 
  (23,6);(19,3)  **\crv{(20,3.5)}, 
 (10,0);(20,12) **@{-}, 
(10,12);(20,0) **@{-}, 
(13,6);(17,6) **@{-}, (13,6.1);(17,6.1) **@{-}, (13,5.9);(17,5.9)
**@{-}, (13,6.2);(17,6.2) **@{-}, (13,5.8);(17,5.8) **@{-}, 
(13,3.1) *{\urcorner}, (17.5,8.9) *{\llcorner}, (19,1.1) *{\lrcorner},  (21,3.5) *{\urcorner},(13,8) *{\ulcorner},
(9,14)*{x \otr y}, (5,8)*{y \utr x}, (8,0)*{x}, (15,0)*{y}, (22,0)*{x}, (25,8)*{y \utr x},
(35,7);(45,7) **@{-} ?>*\dir{>}, 
(35,5);(45,5) **@{-} ?<*\dir{<},
   (63,9.8);(67,9.8)  **\crv{(65,10.3)}, 
  (57,6);(61,9)  **\crv{(60,8.5)}, 
  (73,6);(69,9)  **\crv{(70,8.5)}, 
(70,12);(60,0) **@{-}, 
(70,0);(60,12) **@{-}, 
(63,6);(67,6) **@{-}, (63,6.1);(67,6.1) **@{-}, (63,5.9);(67,5.9)
**@{-},(63,6.2);(67,6.2) **@{-}, (63,5.8);(67,5.8) **@{-}, 
(62,2) *{\urcorner}, (67,8.3) *{\llcorner}, (67.5,3) *{\lrcorner},  (70.9,8) *{\lrcorner},(61.3,10) *{\ulcorner}, 
(-5,6)*{\Gamma_4':},
(60,14)*{x \otr y}, (56,5)*{y \utr x}, (65,3)*{x}, (65,12)*{y}, (77,5)*{y \utr x}, (70,14)*{x \otr y}, 
\endxy}

\vskip.3cm

 \centerline{ \xy (9,2);(13,6) **@{-}, (9,6);(10.5,4.5) **@{-},
(11.5,3.5);(13,2) **@{-}, (17,2);(21,6) **@{-}, (17,6);(21,2)
**@{-}, (13,6);(17,6) **\crv{(15,8)}, (13,2);(17,2) **\crv{(15,0)},
(7,7);(9,6) **\crv{(8,7)}, (7,1);(9,2) **\crv{(8,1)}, (23,7);(21,6)
**\crv{(22,7)}, (23,1);(21,2) **\crv{(22,1)}, 
(17,4);(21,4) **@{-}, (17,4.1);(21,4.1) **@{-}, (17,3.9);(21,3.9)
**@{-}, (17,4.2);(21,4.2) **@{-}, (17,3.8);(21,3.8) **@{-},
(6,7)*{x}, (6,1)*{x}, (27,1)*{x \utr x}, (27,7)*{x \otr x},
(10,3) *{\llcorner},  (12,3) *{\lrcorner}, (21,6) *{\llcorner},(21,2.2) *{\lrcorner},
(35,5);(45,5) **@{-} ?>*\dir{>}, (35,3);(45,3) **@{-} ?<*\dir{<},
(59,2);(63,6) **@{-}, (59,6);(63,2) **@{-}, (67,2);(71,6) **@{-},
(67,6);(68.5,4.5) **@{-}, (69.5,3.5);(71,2) **@{-}, (63,6);(67,6)
**\crv{(65,8)}, (63,2);(67,2) **\crv{(65,0)}, (57,7);(59,6)
**\crv{(58,7)}, (57,1);(59,2) **\crv{(58,1)}, (73,7);(71,6)
**\crv{(72,7)}, (73,1);(71,2) **\crv{(72,1)}, 
(63,4);(59,4) **@{-}, (63,4.1);(59,4.1) **@{-}, (63,3.9);(59,3.9)
**@{-}, (63,4.2);(59,4.2) **@{-}, (63,3.8);(59,3.8) **@{-},
 (56,1)*{x}, (56,7)*{x}, (77,7)*{x \utr x}, (77,1
 )*{x \otr x}, 
(59.5,2.5) *{\llcorner},  (59,6) *{\lrcorner}, (71,6) *{\llcorner},(71,2.2) *{\lrcorner},
 (-5,4)*{\Gamma_5:},
\endxy }

\vskip.3cm

\centerline{ \xy (3,20);(14.8,6.8) **@{-}, (15.2,6.2);(18.8,2.2) **@{-},
(19.5,1.5);(27,-7) **@{-}, (12,20);(3,11) **@{-}, (3,11);(10.8,2.3) **@{-},
(11.3,1.3);(14.4,-2.5) **@{-}, (15.3,-3.8);(18.2,-7) **@{-},
(18,20);(27,11) **@{-}, (12,-7);(27,11) **@{-}, (3,-7);(27,20) **@{-},
(7,17);(7,13) **@{-}, (7.1,17);(7.1,13) **@{-},
(7.4,17);(7.4,13) **@{-}, (7.2,17);(7.2,13) **@{-},
(7.3,17);(7.3,13) **@{-}, (21,15);(25,15) **@{-},
(21,15.1);(25,15.1) **@{-}, (21,15.3);(25,15.3) **@{-},
(21,15.2);(25,15.2) **@{-}, (21,15.4);(25,15.4) **@{-},
 (4.5,18) *{\ulcorner}, (19.5,18) *{\ulcorner}, (5,8.3) *{\ulcorner}, (11,19) *{\llcorner}, (18,0) *{\llcorner}, (26.2,19) *{\llcorner}, (20,11.5) *{\urcorner}, (20.8,0) *{\lrcorner},
 (3,22)*{x}, (27,22)*{y},  (3,-9)*{y}, (11.5,-9)*{y}, (18,-9)*{x}, (27,-9)*{x}, (11,6.8)*{y\otr x},  (20,6.8)*{x\utr y}, (10.5,-4)*{x\utr y},  (20,-4)*{y\otr x},
(35,11);(45,11) **@{-} ?>*\dir{>}, (35,9);(45,9) **@{-} ?<*\dir{<},
(77,20);(65.2,6.8) **@{-}, (64.8,6.2);(61,2) **@{-}, (60.5,1.5);(53,-7)
**@{-}, (68,20);(77,11) **@{-}, (77,11);(69.2,2.2) **@{-},
(68.8,1.8);(65,-2.5) **@{-}, (64.2,-3.2);(61,-7) **@{-},
(62,20);(53,11) **@{-}, (68,-7);(53,11) **@{-}, (77,-7);(53,20)
**@{-},
 (54.5,18) *{\ulcorner}, (69.5,18) *{\ulcorner}, (55,8.4) *{\ulcorner}, (61,19) *{\llcorner}, (67.5,0) *{\llcorner}, (76.2,19) *{\llcorner}, (70,11.5) *{\urcorner}, (70.8,0) *{\lrcorner},
(57,17);(57,13) **@{-}, (57.1,17);(57.1,13) **@{-},
(57.2,17);(57.2,13) **@{-}, (57.3,17);(57.3,13) **@{-},
(57.4,17);(57.4,13) **@{-}, (71,15);(75,15) **@{-},
(71,15.1);(75,15.1) **@{-}, (71,15.3);(75,15.3) **@{-},
(71,15.2);(75,15.2) **@{-}, (71,15.4);(75,15.4) **@{-},
(-5,10)*{\Gamma_{8}:},
 (53,22)*{x}, (77,22)*{y},  (53,-9)*{y}, (61.5,-9)*{y},
 (68,-9)*{x}, (77,-9)*{x}, (60.5,6.8)*{y\utr x},  (69.5,6.8)*{x\otr y}, (59.6,-4)*{x\otr y},  (70,-4)*{y\utr x},
\endxy}

The biquandle axioms are chosen so that for a given biquandle coloring of 
a diagram
on one side of a move (Reidemeister move in the case of classical links,
Yoshikawa move in the case of surface-links) there is a unique biquandle
coloring of the diagram on the other side of the move. Hence by construction
we have the following standard result:

\begin{theorem}
Let $X$ be a finite biquandle and $L$, $L'$ marked graph diagrams of 
ambient isotopic surface-links. Then the number of $X$-colorings of $L$ 
and the number of $X$-colorings of $L'$ are equal. 
\end{theorem}

\begin{definition}
Let $X$ be a finite biquandle. The number of $X$-colorings of a surface-link
represented by a marked graph diagram $L$ is called the \textit{biquandle
counting invariant} of $L$ with respect to $X$, denoted by
$\Phi_X^{\mathbb{Z}}(L)$.
\end{definition}

\begin{remark}
We can also define $\Phi_X^{\mathbb{Z}}(L)$ as the cardinality of the set of
\textit{biquandle homomorphisms} $f:\mathcal{B}(L)\to X$, i.e.,
maps satisfying
\[f(x\utr y)=f(x)\utr f(y)\quad \mathrm{and}\quad 
f(x\otr y)=f(x)\otr f(y)
\]
for all $x,y\in \mathcal{B}(L).$
\end{remark}

\section{Biquandle Module Enhancements}\label{C}

We will now adapt an idea from previous work (see for example \cite{BN}) to 
enhance the biquandle counting invariant for surface-links.

\begin{definition}
Let $X$ be a finite biquandle and let $R$ be a commutative ring
with identity. A \textit{biquandle module} over $X$ with coefficients
in $R$ is an assignment of units $t_{x,y},r_{x,y}\in R^{\times}$ and 
elements $s_{x,y}\in R$ satisfying for all $x,y,z\in X$
\[\begin{array}{rcll}
t_{x,x}+s_{x,x} & = & r_{x,x}, & (i.i) \\
r_{y\otr x,z\otr x} r_{x,z} & = & r_{x\utr y, z\otr y}r_{y,z}, & (iii.i) \\
r_{x\utr z,y\utr z} t_{y,z} & = & t_{y\otr x, z\otr x}r_{x,y}, & (iii.ii) \\
r_{x\utr z,y\utr z} s_{y,z} & = & s_{y\otr x, z\otr x}r_{x,z}, & (iii.iii) \\
t_{x\utr z,y\utr z} t_{x,z} & = & t_{x\utr y, z\otr y}t_{x,y}, & (iii.iv) \\
s_{x\utr z,y\utr z} t_{y,z} & = & t_{x\utr y, z\otr y}s_{x,y}, & (iii.v) \\
t_{x\utr z,y\utr z} s_{x,z} +s_{x\utr z, y\utr z}s_{y,z} & = & s_{x\utr y, z\otr y}r_{y,z}.  & (iii.vi)
\end{array}\]
\end{definition}

We can specify an $X$-module with a triple $[t,s,r]$ of matrices 
$t,s,r\in M_n(R)$ whose row $x$ column $y$ entries are $t_{x,y}, s_{x,y}$ 
and $r_{x,y}$ respectively.

\begin{example}\label{ex1}
Let $X$ be the biquandle structure on the set $X=\{1,2\}$ specified by the 
operation tables:
\[\begin{array}{r|rr}
\utr & 1 & 2 \\ \hline
1 & 2 & 2 \\
2 & 1 & 1
\end{array}
\quad 
\begin{array}{r|rr}
\otr & 1 & 2 \\ \hline
1 & 2 & 2 \\
2 & 1 & 1
\end{array}
\]
Then our \texttt{python} computations reveal $X$-module structures over 
$\mathbb{Z}_5$ including
\[\left[
\left[\begin{array}{rr} 
2 & 3 \\
4 & 1
\end{array}\right],
\left[\begin{array}{rr} 
2 & 0\\
0 & 1
\end{array}\right],
\left[\begin{array}{rr} 
4 & 4 \\
3 & 2
\end{array}
\right]\right].\]
Then for instance we have $t_{1,1}=2$, $s_{1,1}=2$ and $r_{1,1}=4$, and 
for $x=1,y=2, z=2$ we have for axiom (iii.vi)
\[t_{1\utr 2,2\utr 2}s_{1,2}+s_{1\utr 2,2\utr 2}s_{2,2}=
t_{2,1}s_{1,2}+s_{2,1}s_{2,2}=4(0)+0(1)=0=0(2)
=s_{2,1}r_{2,2}
=s_{1\utr 2,2\otr 2}r_{2,2}
\]
and for axiom (iii.iii)
\[r_{1\utr 2,2\utr 2}s_{2,2}
=r_{2,1}s_{2,2}
=3(1)=3
=(2)(4)
=s_{1,1}r_{1,2}
=s_{2\otr 1,2\otr 1}r_{1,2},
\]
etc.
\end{example}

The biquandle module axioms are motivated by the idea of enhancing biquandle
colorings of oriented link diagrams with secondary ``bead colorings'' where 
the beads satisfy equations similar to the Alexander biquandle relations
but with coefficients which depend on the biquandle colors at a crossing:
\[\includegraphics{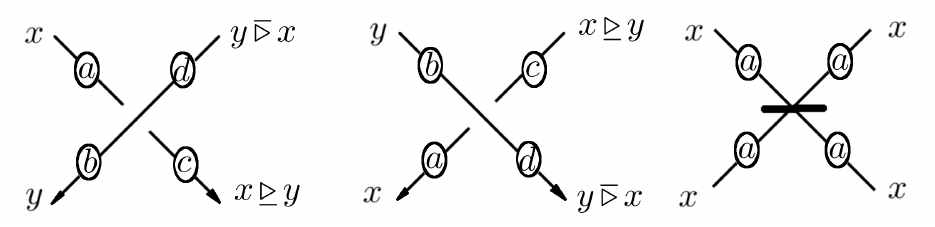}\]
\[\begin{array}{rcl}
c & = & t_{x,y} a+s_{x,y} b, \\
d & = & r_{x,y} b.
\end{array}\]

The condition (i.i), i.e.,
\[t_{x,x}+s_{x,x}=r_{x,x},\]
 is required by Reidemeister move I:
\[\includegraphics{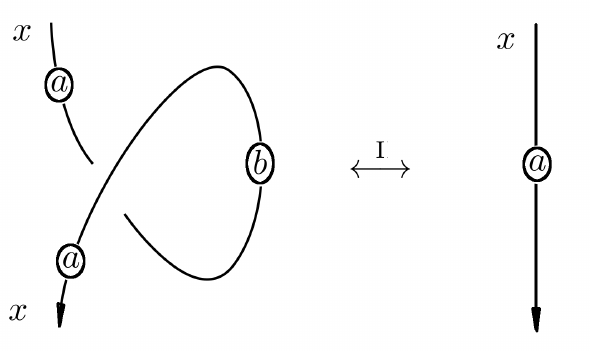}\]
which is satisfied provided that $b=(t_{x,x}+s_{x,x})a=r_{x,x}a$.

The requirement that $t_{x,y}$ and $r_{x,y}$ are invertible implies that
the pair $(c,d)$ determines the pair $(a,b)$, the pair $(a,d)$ determines
the pair $(b,c)$ and the pair $(b,c)$ determines the pair $(a,d)$. This
suffices to guarantee that a bead coloring on one side of a Reidemeister II
move corresponds to a unique bead coloring on the other side of the move.
\[\includegraphics{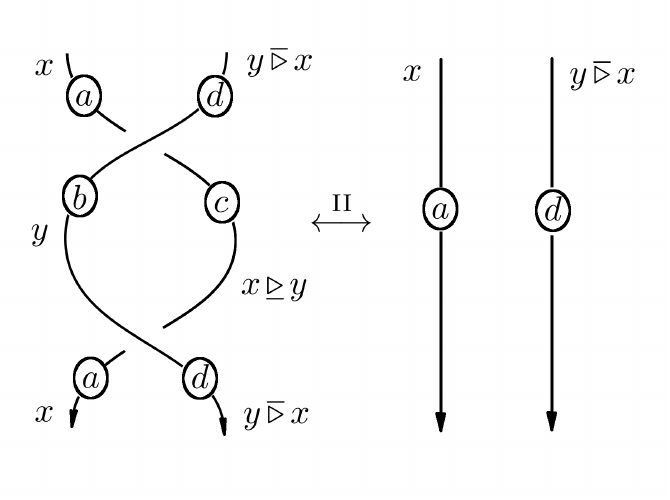}\quad\raisebox{0.3in}{
\includegraphics{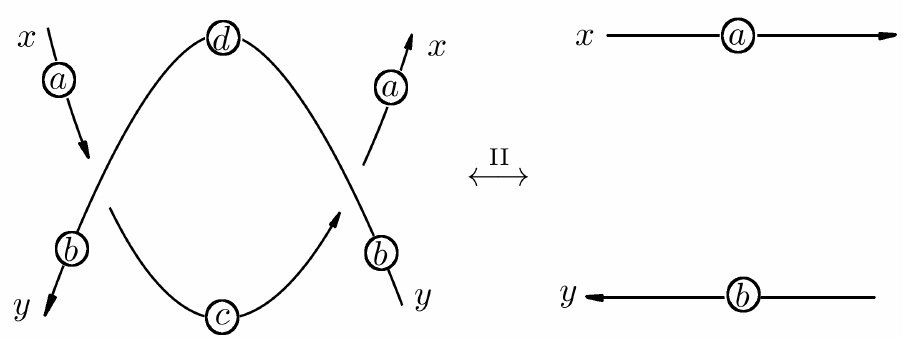}}
\]

The Reidemeister III move yields the conditions (iii.i) through (iii.vi):
\begin{eqnarray*}
f & = & r_{x\utr y, z\otr y} j = r_{x\utr y, z\otr y} r_{y,z} c \\
 & = & r_{y\otr x,z\otr x} k = r_{y\otr x,z\otr x} r_{x,z} c, \\
h & = & t_{y\otr x,z\otr x} d+s_{y\otr x,z\otr x} k 
= t_{y\otr x,z\otr x} r_{x,y}b+s_{y\otr x,z\otr x} r_{x,z} c \\
& = & r_{x\utr z,y\utr z} e = r_{x\utr z,y\utr z} t_{y,z}b+r_{x\utr z,y\utr z} s_{y,z}c, \\ 
g & = & t_{x\utr y,z\otr y}i+s_{x\utr y,z\otr y}j 
=t_{x\utr y,z\otr y}t_{x,y}a+t_{x\utr y,z\otr y}s_{x,y}b+s_{x\utr y,z\otr y}r_{y,z}c\\
 & = & t_{x\utr z,y\utr z}l+s_{x\utr z,y\utr z}e
=t_{x\utr z,y\utr z}t_{x,z}a+s_{x\utr z,y\utr z}t_{y,z}b+
(t_{x\utr z,y\utr z}s_{x,z}+s_{x\utr z,y\utr z}s_{y,z})c.
\end{eqnarray*}
\[\includegraphics{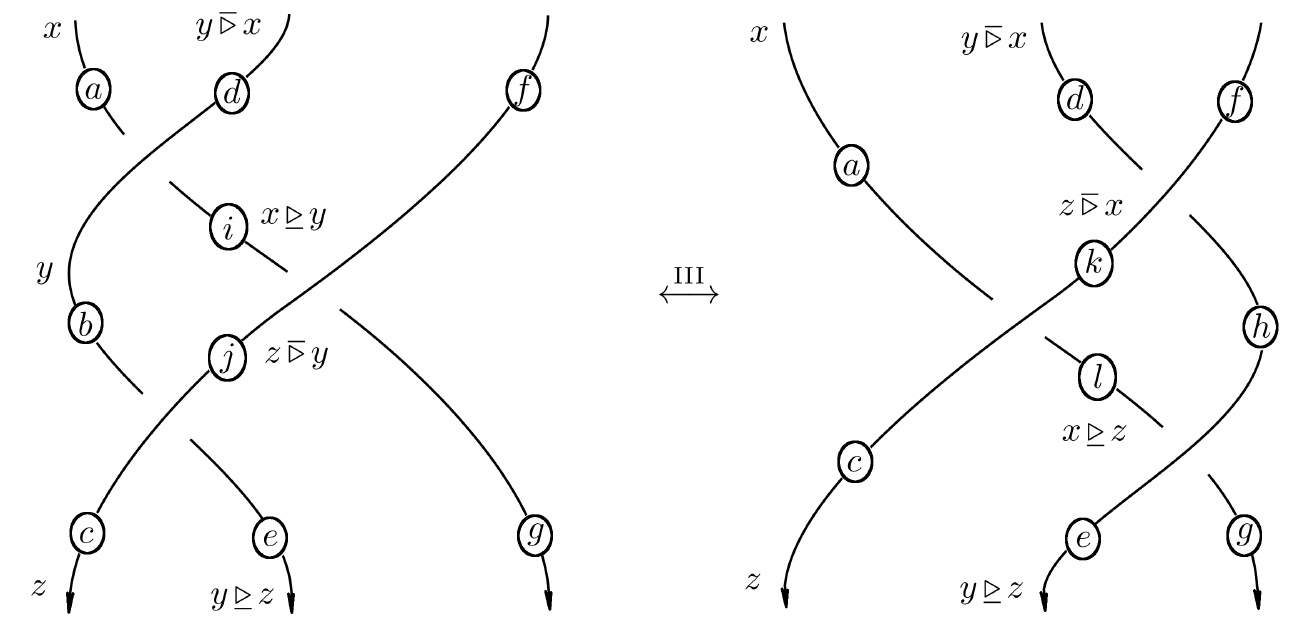}\]

The moves $\Gamma_4,\Gamma_4',\Gamma_6$ and $\Gamma_7$ do not impose additional
conditions on the bead colorings.
\[\includegraphics{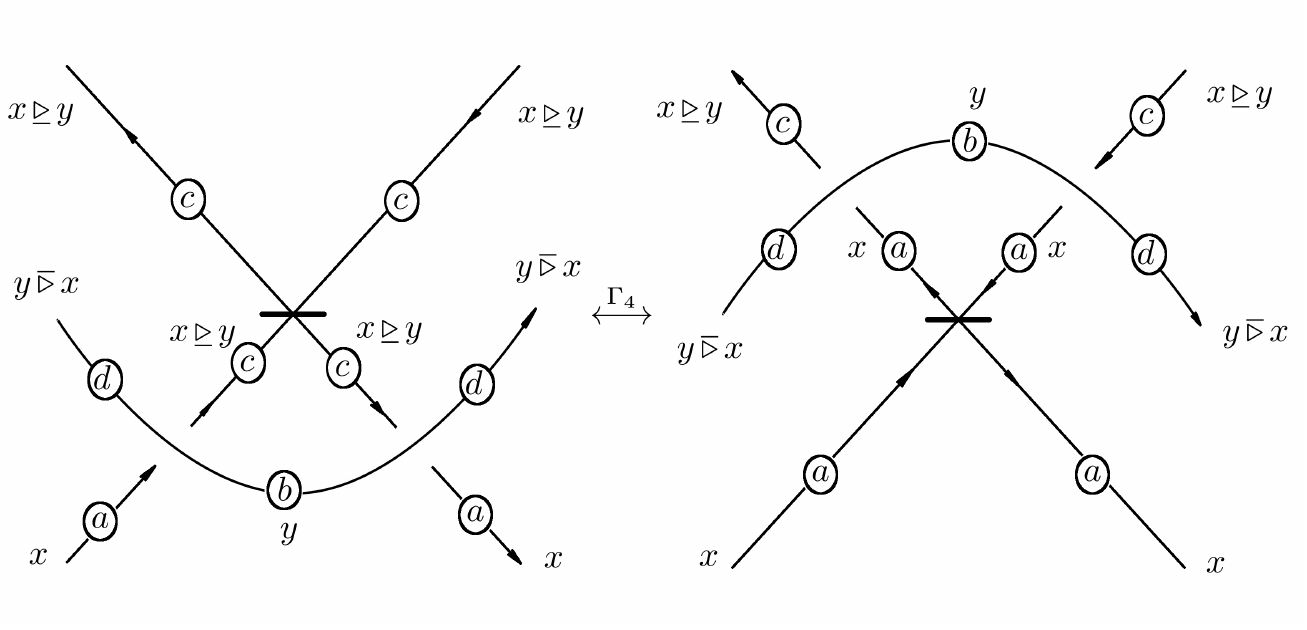}\]
\[\includegraphics{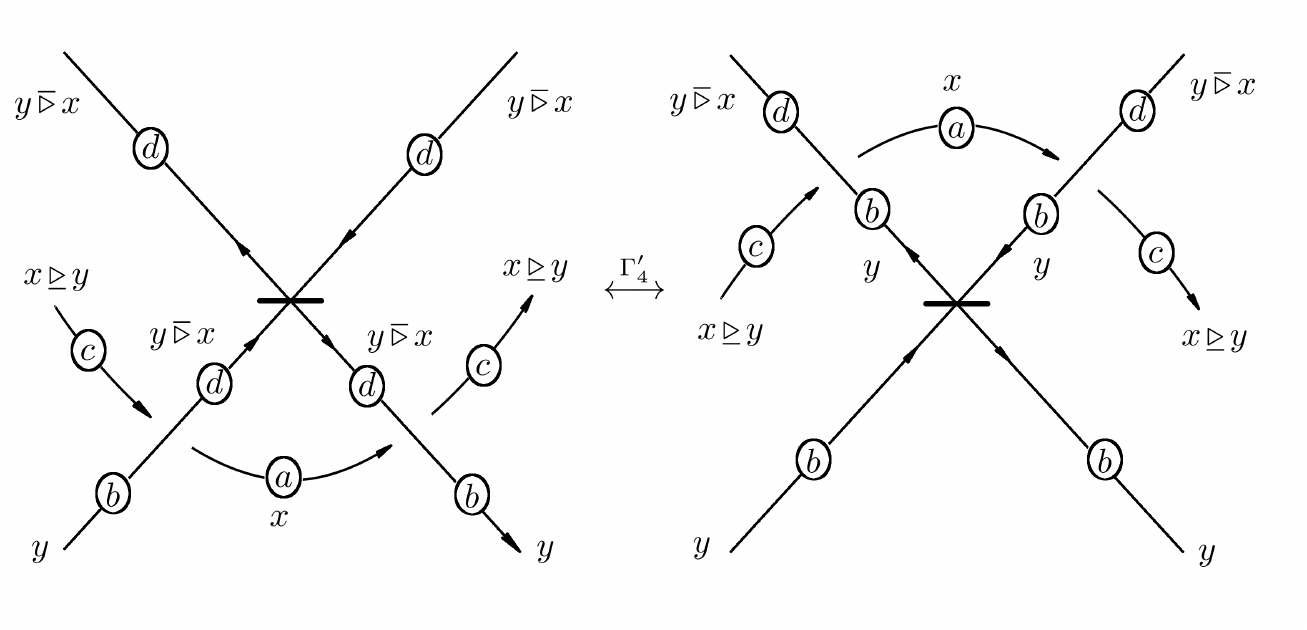}\]
\[\includegraphics{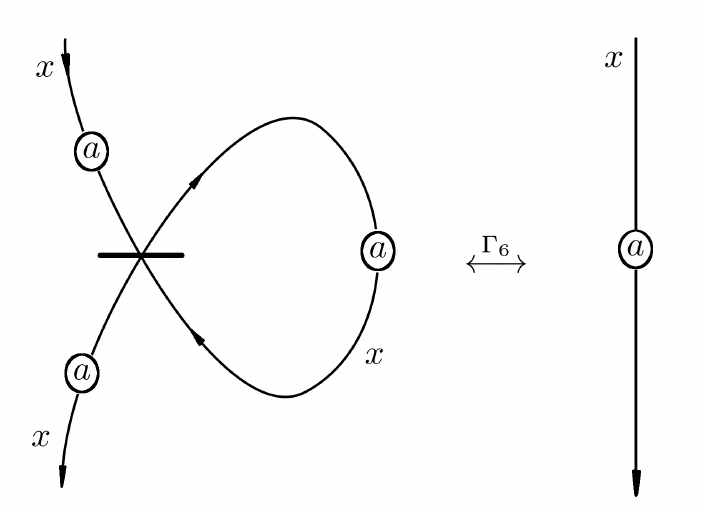}\]
\[\includegraphics{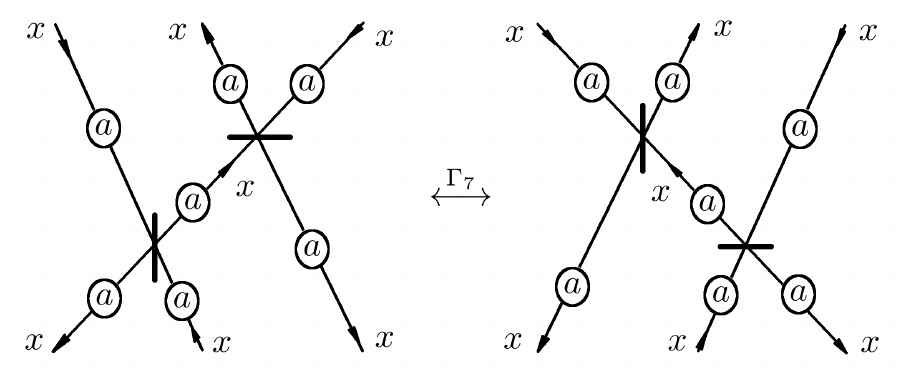}\]

The condition (i.i) also ensures that the space of bead colorings
is unchanged by Yoshikawa $\Gamma_5$-moves:
\[\includegraphics{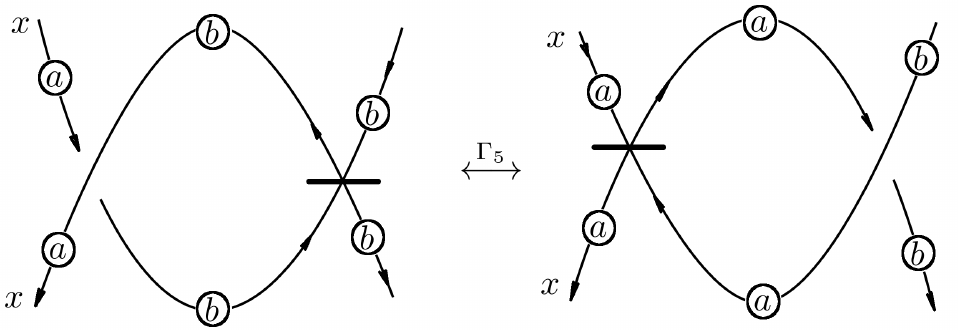}\]

Finally, the Yoshikawa $\Gamma_8$-move does not require any additional 
conditions; the beads 
$c=t_{x,y}^{-1}a-t_{x,y}^{-1}s_{x,y}b$, 
$d=r_{x,y}^{-1}b$, 
$e=t_{y,x}^{-1}b-t_{y,x}^{-1}s_{y,x}a$ 
and $f=r_{y,x}^{-1}a$ 
are uniquely determined by the biquandle colors $x,y\in X$ 
and beads $a,b\in R$:
\[\includegraphics{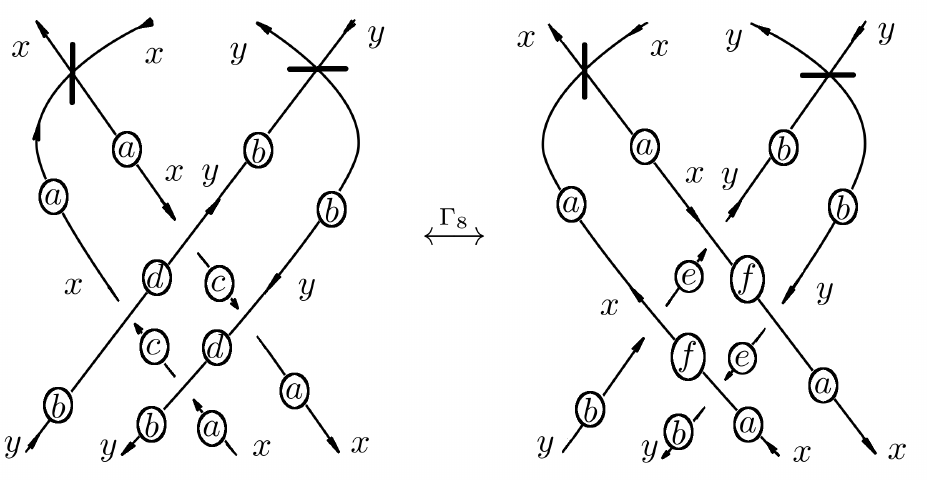}\]

Thus by construction we have the following:

\begin{theorem}
Let $X$ be a finite biquandle, $R$ a commutative ring, $m=[t,s,r]$ an
$X$-module over $R$, and $L$ an oriented marked graph diagram. Then
for each $X$-coloring $L_f$ of $L$, the number $\phi_m(L_f)$ of bead colorings
by $m$ is unchanged by $X$-colored Yoshikawa moves.
\end{theorem}

\begin{corollary}
Let $X$ be a finite biquandle, $R$ a commutative ring, $[t,s,r]$ an
$X$-module over $R$, and $L$ an oriented marked graph diagram. Then the
multiset 
\[\Phi_X^{M,m}(L)=\{|\phi_m(L_f)|\ f\in\mathrm{Hom}(\mathcal{B}(L),X)\}\]
and polynomial
\[\Phi_X^{m}(L)=\sum_{f\in\mathrm{Hom}(\mathcal{B}(L),X)}u^{|\phi_m(L_f)|} \]
are invariants of surface-links called the \textit{biquandle module multiset}
and \textit{biquandle module polynomial} respectively.
\end{corollary}

\begin{example}
Let us illustrate the computation of the invariant for the surface-link 
$L=6_1^{0,1}$ below using the biquandle $X$ and module $m$ from Example 
\ref{ex1}. This module data gives us bead coloring rules 
\[\includegraphics{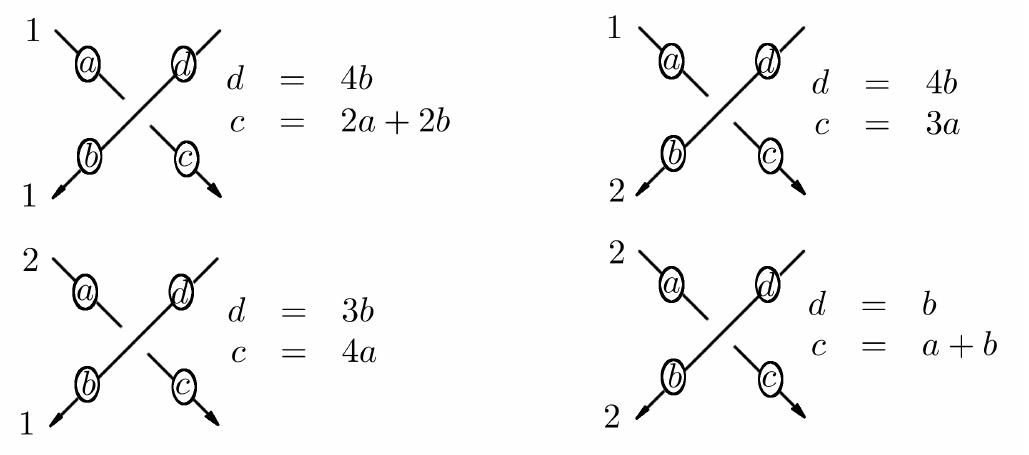}\]
where the beads are elements of $\mathbb{Z}_5$.
There are four $X$-colorings of $L$ as shown. 
\[\includegraphics{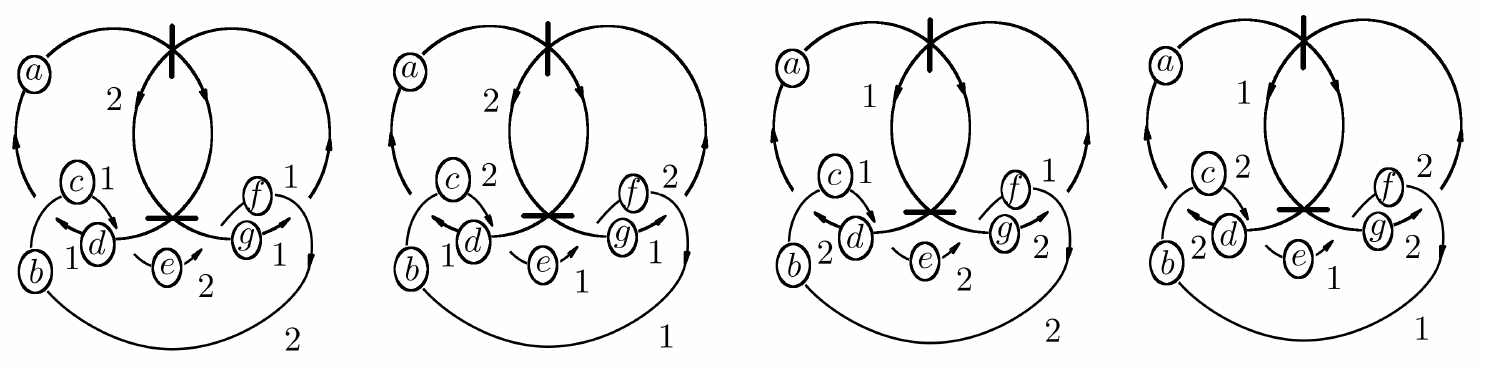}\]
The first $X$-coloring has the following system of bead coloring equations, 
yielding the following coloring matrix: 
\[\begin{array}{rcl}
t_{11}d+ s_{11}c & = & a\\
r_{11}c & = & b\\
t_{11}c+ s_{11}d & = & e\\
r_{11}d & = & a \\
t_{21}e +s_{21}g &=& f\\
r_{21}g & = & a \\
s_{12} b+t_{12} g & = & a \\
r_{12}b & = & f
\end{array}
\leftrightarrow
\left[\begin{array}{rrrrrrr}
-1 & 0 & s_{11} & t_{11} & 0 & 0 & 0 \\
0 & -1 & r_{11} & 0 & 0 & 0 & 0 \\
0 & 0 & t_{11} & s_{11} & -1 & 0 & 0 \\
-1 & 0 & 0 & r_{11} & 0 & 0 & 0 \\
0 & 0 & 0 & 0 & t_{21} & -1 & s_{21} \\
-1 & 0 & 0 & 0 & 0 & 0 & r_{21} \\
-1 & s_{12} & 0 & 0 & 0 & 0 & t_{12} \\
0 & r_{12} & 0 & 0 & 0 & -1 & 0 \\
\end{array}\right].
\]
Then substituting the values from our chosen $X$-module we obtain
matrix over $\mathbb{Z}_5$
\[
\left[\begin{array}{rrrrrrr}
4 & 0 & 2 & 2 & 0 & 0 & 0 \\
0 & 4 & 4 & 0 & 0 & 0 & 0 \\
0 & 0 & 2 & 2 & 4 & 0 & 0 \\
4 & 0 & 0 & 4 & 0 & 0 & 0 \\
0 & 0 & 0 & 0 & 4 & 4 & 0 \\
4 & 0 & 0 & 0 & 0 & 0 & 3 \\
4 & 0 & 0 & 0 & 0 & 0 & 3 \\
0 & 4 & 0 & 0 & 0 & 4 & 0 \\
\end{array}\right]
\leftrightarrow
\left[\begin{array}{rrrrrrr}
1 & 0 & 0 & 0 & 0 & 0 & 2 \\
0 & 1 & 0 & 0 & 0 & 0 & 2 \\
0 & 0 & 1 & 0 & 0 & 0 & 3 \\
0 & 0 & 0 & 1 & 0 & 0 & 3 \\
0 & 0 & 0 & 0 & 1 & 0 & 2 \\
0 & 0 & 0 & 0 & 0 & 1 & 3 \\
0 & 0 & 0 & 0 & 0 & 0 & 0 \\
0 & 0 & 0 & 0 & 0 & 0 & 0 \\
\end{array}\right]\]
and this $X$-coloring has a one dimensional space of bead colorings
\[\includegraphics{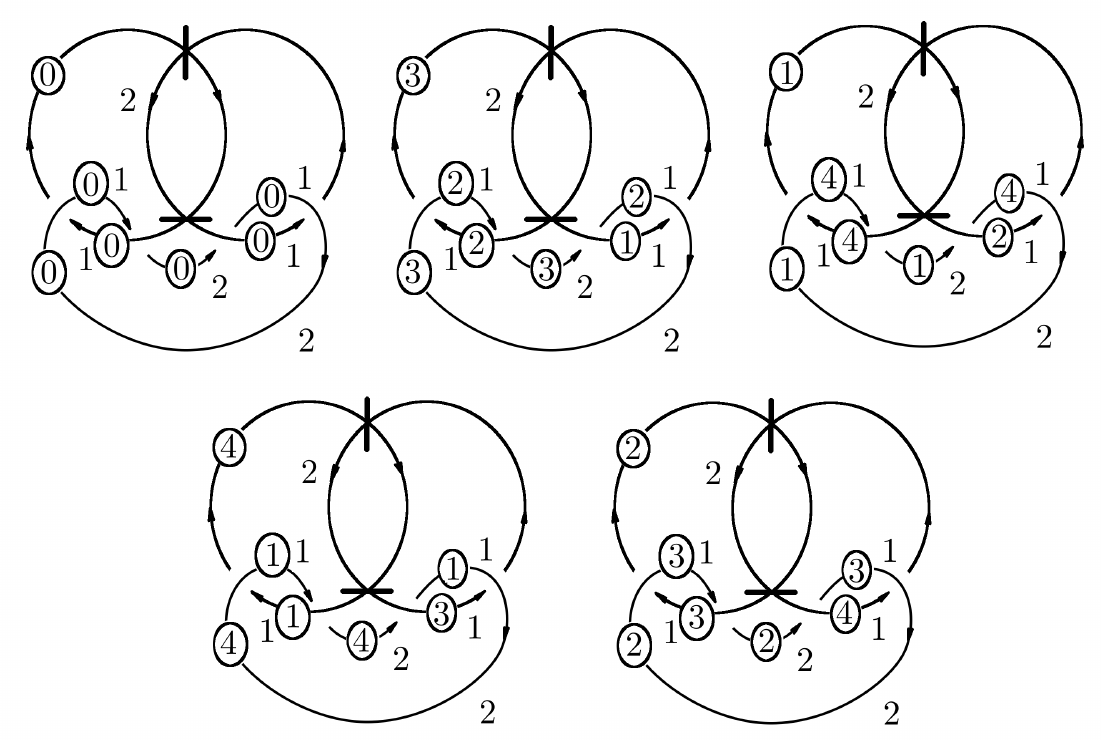}\]
contributing $u^5$ to the invariant. Computing the other
three, we obtain the invariant value $2u^5+2u^{25}$. We observe that this contains
more information than the biquandle counting invariant 
$\Phi_X^{\mathbb{Z}}(6_1^{0,1})=4$ alone, since it separates
the colorings into two sets: two with 25 bead colorings and two with 5.
In particular, this invariant distinguishes this surface-link from the unlink
of a torus and sphere which has the invariant value $4u^5$.
\end{example}

\begin{example}
Let $X$ be the biquandle 
\[\begin{array}{r|rr}
\utr & 1 & 2 \\ \hline
1 & 2 & 2 \\
2 & 1 & 1
\end{array}
\quad 
\begin{array}{r|rr}
\otr & 1 & 2 \\ \hline
1 & 2 & 2 \\
2 & 1 & 1
\end{array}
\]
from Example \ref{ex1} and $m$ the biquandle module
over $\mathbb{Z}_5$ given by the matrices
\[\left[
\left[\begin{array}{rr} 3 & 4 \\ 4 & 3 \end{array}\right],
\left[\begin{array}{rr} 0 & 0 \\ 0 & 0\end{array}\right],
\left[\begin{array}{rr} 3 & 1 \\ 1 & 3\end{array}\right]\right].
\]
Our \texttt{python} computations give the following $\Phi_X^{m}$ values
for orientable surface-links of small $ch$-index with orientations as shown.
\[\includegraphics{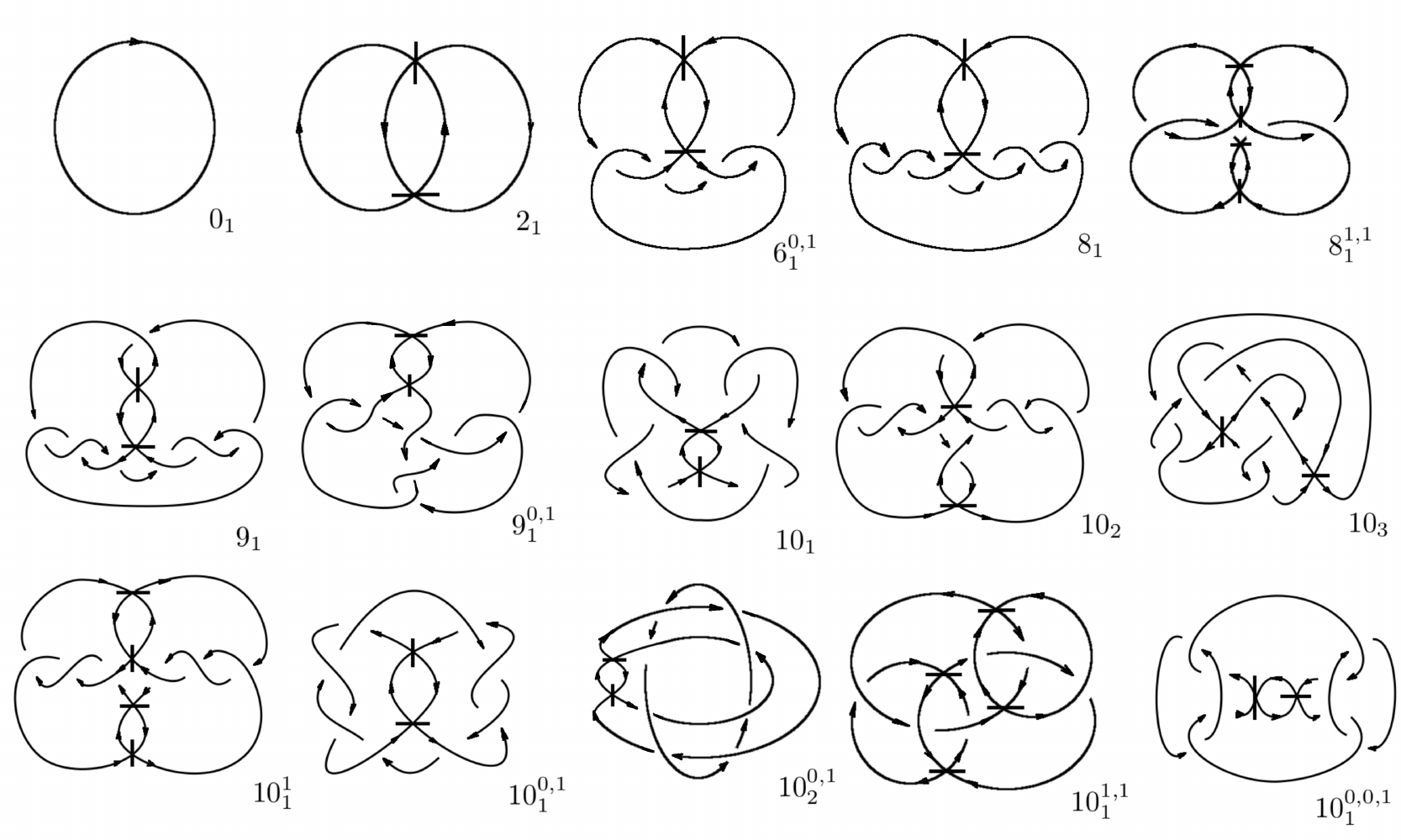}\]
The results are in the table:
\[
\begin{array}{r|l}
\Phi_X^m(L) & L \\\hline
2u^5 & 8_1, 9_1,10_1, 10_2, 10_3 \\
2u^5+2u^{25} & 6^{0,1}_1,10^{0,1}_2 \\
2u+2u^{25} & 8^{1,1}_1, 10_1^{1,1} \\
4u^{25} & 9_1^{0,1}, 10^{0,1}_1 \\
4u^{25}+4u^{125} & 10_1^{0,0,1}.
\end{array}
\]
We can observe that this particular pair of biquandle and module do not 
distinguish the
surface-knots in this small sample, but are effective at distinguishing the 
surface-links from each other.
\end{example}

\begin{example}
For our final example we computed $\Phi_X^{m}$ for the orientable surface-links
of small $ch$-index with respect to the $X$-modules over $\mathbb{Z}_3$ 
given by the matrices
\[m_1=
\left[
\left[\begin{array}{rrr} 2 & 2 & 2 \\ 2 & 2 & 2 \\ 1 & 1 & 1\end{array}\right],
\left[\begin{array}{rrr} 0 & 0 & 0 \\ 0 & 0 & 0 \\ 0 & 0 & 0\end{array}\right],
\left[\begin{array}{rrr} 2 & 2 & 1 \\ 1 & 2 & 2 \\ 1 & 1 & 1\end{array}\right]\right]
\]
and
\[m_2=
\left[
\left[\begin{array}{rrr} 1 & 1 & 1 \\ 1 & 1 & 1 \\ 2 & 2 & 2\end{array}\right],
\left[\begin{array}{rrr} 1 & 2 & 2 \\ 2 & 1 & 1 \\ 1 & 2 & 2\end{array}\right],
\left[\begin{array}{rrr} 2 & 1 & 2 \\ 2 & 2 & 1 \\ 1 & 1 & 1\end{array}\right]\right]
\]
for the biquandle 
\[
\begin{array}{r|rrr}
\utr & 1 & 2 & 3 \\ \hline
1 & 2 & 2 & 2  \\
2 & 1 & 1 & 1 \\
3 & 3 & 3 & 3 
\end{array}
\quad
\begin{array}{r|rrr}
\otr & 1 & 2 & 3 \\ \hline
1 & 2 & 3 & 1  \\
2 & 3 & 1 & 2 \\
3 & 1 & 2 & 3 
\end{array}
\]
The results are in the table.
\[
\begin{array}{lll}
L& \Phi_X^{m_1}(L) & \Phi_X^{m_2}(L)  \\ \hline
2_1 & 3u^3 & 3u^3 \\ 
6_1^{0,1} & 3u^9 & 3u^3 \\
8_1 & 9u^3 & 9u^9 \\
8_1^{1,1} & 3u^9  & 3u^3 \\
9_1 & 9u^3 & 9u^9 \\
9_1^{0,1} & 3u^9 & 3u^3 \\
10_1 & 3u^3 & 3u^3 \\
10_2 & 9u^3 & 9u^9 \\
10_3 & 3u^3 & 3u^3 \\
10^1_1 & 9u^3 & 9u^9 \\
10_1^{0,1} & 3u^9 & 3u^3 \\
10_2^{0,1} & 3u^9 & 3u^9 \\
10_1^{1,1} & 3u^9 & 3u^3 \\
10_1^{0,0,1} & 9u^{27} & 9u^9 
\end{array}
\]
\end{example}

\section{Questions}\label{Q}

We end with some questions and directions for future research. 

Faster methods for finding biquandle modules would be desirable; our current
approach fills in entries in the $[t,s,r]$ matrix using the module conditions 
and works well enough for small biquandles and small rings, but other methods 
will be necessary for finding biquandle modules over larger finite and 
infinite rings.

As in \cite{CN}, biquandle modules over polynomial rings should be of interest.
Such a module effectively defines Alexander invariants for biquandle-colored
links. This suggests natural questions such as:
\begin{itemize}
\item Which, if any, skein relations are satisfied by various biquandle 
modules?
\item What kinds of categorifications can be defined for these invariants?
\item What additional enhancements can be defined in the case of biquandle module
invariants of surface-links?
\end{itemize}
et cetera.

\bibliography{yj-sn}{}

\begin{thebibliography}{10}

\bibitem{AG}
N.~Andruskiewitsch and M.~Gra\~{n}a.
\newblock From racks to pointed {H}opf algebras.
\newblock {\em Adv. Math.}, 178(2):177--243, 2003.

\bibitem{BN}
R.~Bauernschmidt and S.~Nelson.
\newblock Birack modules and their link invariants.
\newblock {\em Commun. Contemp. Math.}, 15(3):1350006, 13, 2013.

\bibitem{BKLNS}
J.~Blankstein, S.~Kim, C.~Lepel, S.~Nelson, and N.~Sanderson.
\newblock Virtual shadow modules and their link invariants.
\newblock {\em Internat. J. Math.}, 23(9):1250096, 22, 2012.

\bibitem{CEGS}
J.~S. Carter, M.~Elhamdadi, M.~Gra\~{n}a, and M.~Saito.
\newblock Cocycle knot invariants from quandle modules and generalized quandle
  homology.
\newblock {\em Osaka J. Math.}, 42(3):499--541, 2005.

\bibitem{CN}
E.~Cody and S.~Nelson.
\newblock Polynomial birack modules.
\newblock {\em Topology Appl.}, 173:285--293, 2014.

\bibitem{GN}
M.~Grier and S.~Nelson.
\newblock Kei modules and unoriented link invariants.
\newblock {\em Homology Homotopy Appl.}, 16(1):167--177, 2014.

\bibitem{HHNYZ}
A.~Haas, G.~Heckel, S.~Nelson, J.~Yuen, and Q.~Zhang.
\newblock Rack module enhancements of counting invariants.
\newblock {\em Osaka J. Math.}, 49(2):471--488, 2012.

\bibitem{KJL}
J.~Kim, Y.~Joung, and S.~Y. Lee.
\newblock On the {A}lexander biquandles of oriented surface-links via marked
  graph diagrams.
\newblock {\em J. Knot Theory Ramifications}, 23(7):1460007, 26, 2014.

\bibitem{KJL2}
J.~Kim, Y.~Joung, and S.~Y. Lee.
\newblock On generating sets of {Y}oshikawa moves for marked graph diagrams of
  surface-links.
\newblock {\em J. Knot Theory Ramifications}, 24(4):1550018, 21, 2015.

\bibitem{NP}
S.~Nelson and K.~Pelland.
\newblock Birack shadow modules and their link invariants.
\newblock {\em J. Knot Theory Ramifications}, 22(10):1350056, 12, 2013.

\bibitem{Y}
K.~Yoshikawa.
\newblock An enumeration of surfaces in four-space.
\newblock {\em Osaka J. Math.}, 31(3):497--522, 1994.

\end{thebibliography}
\bibliographystyle{abbrv}

\bigskip

\noindent
\textsc{Department of Mathematics \\
Michigan State University \\
C120, Wells Hall \\
619 Red Cedar Rd \\
East Lansing, MI 48824
}

\medskip

\noindent
\textsc{Department of Mathematical Sciences \\
Claremont McKenna College \\
850 Columbia Ave. \\
Claremont, CA 91711}

\end{document}